\documentclass[11pt,a4paper]{article}
\usepackage{fouriernc}
\usepackage{tabularx}
\usepackage{booktabs}
\usepackage{makecell}
\usepackage{amsmath, amssymb}
\usepackage{graphicx}
\usepackage{xcolor,ulem}
\usepackage{siunitx}
\usepackage[export]{adjustbox}
\usepackage{subfig}
\usepackage[utf8]{inputenc}
\usepackage[T1]{fontenc}
\usepackage{enumitem}
\usepackage{color, soul}
\usepackage{authblk}
\usepackage{marvosym}
\usepackage{wasysym}
\usepackage{multirow}
\usepackage{url}
\usepackage[sort&compress,round,comma,authoryear]{natbib}
\usepackage{hyperref}
\hypersetup{colorlinks=true,allcolors=blue}
\usepackage{vmargin}
\setmarginsrb{2cm}{2cm}{2cm}{2cm}{0mm}{0mm}{0mm}{10mm}
\usepackage{tikz}
\usepackage{pgfplots, pgfplotstable}
\pgfplotsset{compat=1.7}

\title{\vspace*{-1.5cm} \bfseries Teaching mathematical modeling for sustainability:\\ Enhancing interdisciplinary skills in students}
\author{N. Karjanto\thanks{\Letter: \url{natanael@skku.edu}\href{https://orcid.org/0000-0002-6859-447X}{\includegraphics[scale=0.08]{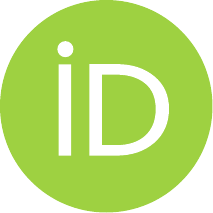}}}}
\affil{Department of Mathematics, University College, Natural Science Campus, Sungkyunkwan University\protect \\ 2066~Seobu-ro, Jangan-gu, Suwon~16419, Gyeonggi-do, Republic of Korea}

\date{\vspace*{-0.5cm} \footnotesize Last updated \today}

\begin{document}
\maketitle
\vspace*{-0.5cm}
\begin{abstract}
\noindent
We developed a pilot course focused on mathematical modeling within the tertiary education framework, with a distinct emphasis on sustainability and sustainable development. While an applicable textbook exists for this liberal arts course, it is noticeable that numerous examples within it are not directly aligned with sustainability concerns. To address this gap, our study strategically integrated the teaching and learning of modeling by carefully selecting classroom examples that closely align with the context of sustainability. By employing an innovative and adaptable approach to the course content delivery, we fostered interdisciplinary collaboration among students, improved their comprehension, and enhanced their interdisciplinary skills.\\

\noindent
Keywords: mathematical modeling, sustainability, liberal arts course, teaching and learning, tertiary education, sustainable development.
\end{abstract}

\section{Introduction}
This paper delves into the realm of sustainable learning and its alignment with the United Nations Sustainable Development Goal 4 (SDG4) in education. Sustainable learning, which encompasses environmental, economic, and social factors, equips learners with tools to navigate a changing world. By integrating sustainability principles, individuals become catalysts for change, fostering a more sustainable society. The integration of sustainable learning practices also aligns with SDG4's goal of inclusive and equitable quality education. 

We emphasize the crucial role of mathematics education in achieving SDG4. Mathematics equips individuals with critical thinking and problem-solving skills vital for personal growth and societal contributions. Given the extensive practical applications of mathematics across diverse fields and the noticeable dearth of studies focusing on the intersection of mathematical modeling and sustainability, we put forth the concept of a liberal arts course in mathematical modeling for sustainability. Our intention is twofold: to nurture sustainable learning in students while addressing sustainability-related concerns through quantitative analysis using mathematical models. This initiative not only aims to foster an in-depth comprehension of mathematical modeling but also endeavors to establish a significant connection between mathematical concepts and sustainability, thereby promoting a comprehensive perspective within mathematics education~\citep{keyley2016a,kaiser2017the}.  	

Although a specific textbook is available and can be adopted for teaching and learning this liberal arts course, it seems that many examples are not (directly) related to issues in sustainability. In this study, we leveraged the teaching and learning of modeling through selected classroom examples that closely align with the context of sustainability. Hence, we are addressing this research question: What innovative approaches can be employed to optimize the teaching of mathematical modeling for sustainability and significantly improve students' comprehension and application of this vital interdisciplinary skill?

\section{Methodology}

We performed an extensive review of mathematical modeling textbooks, eventually selecting one as the primary textbook for our course. We meticulously analyzed this chosen textbook to determine the alignment of its examples and exercises with sustainability principles. Additionally, we explored supplementary textbooks for references, identifying relevant sustainability-related examples to incorporate into our classroom discussions. Crafting suitable examples ourselves also proved to be essential, which contributed to the demanding and time-intensive nature of teaching this pilot course.

\section{Results}

Regarding the potential textbooks for the course, we narrowed down from over 30 into two, that is, \cite{roe2018mathematics} and \cite{hersh2006mathematical}, and chose the former for our class. However, we discovered that not all examples are related to sustainability, and thus further effort is needed to deliver good quality teaching and learning.

Furthermore, we identified several useful examples that can be utilized for a course in mathematical modeling for sustainability. For example, \cite{caldwell2004mathematical} contained one exercise problem on lake pollution. \cite{meerschaert2013mathematical} provided one example of whether competitive dynamics can coexist sustainably in an unmanaged forest area, and another one revolves around the potential extinction of the blue whale species due to competition with the fin whale, taking into account their growth rates and the available environmental resources.

Although \cite{hughes2018applied} is not a textbook on mathematical modeling, it contains several sustainability examples, such as marine harvesting, lake pollution, predator-prey model, and modeling the spread of a disease. \cite{larson2017calculus} offers an example of modeling a chemical mixture, which can be tailored to issues in sustainability, such as the presence of chromium, a hard metal commonly found in the freshwater supply. The following example illustrates such a scenario, with instructor-adjusted intervention:\\

\begin{minipage}[c]{0.93\textwidth}
{\small
Chromium is a toxic heavy metal commonly found in freshwater supplies, and its presence can be harmful to human health and the environment. In an effort to reduce the amount of chromium in a reservoir, a solar-powered water pump is used to transfer water from a well to the reservoir. Initially, the reservoir contains 8000 liters of water with 2~g of chromium dissolved. The pump draws water from the well, which has a concentration of 30~micrograms of chromium per liter, at a rate of 10~liters/minute. The water in the reservoir is then used to irrigate a field of crops, specifically cauliflower, which requires 10~liters of water per minute to grow. After one hour, how much chromium accumulates in the reservoir? What is the steady- state amount of chromium in the reservoir?
\par}
\end{minipage}

\section{Discussion}

Despite the availability of a textbook, designing a pilot liberal arts course on mathematical modeling for sustainability turned out to be challenging. It implies that there are many open problems for both theoreticians and practitioners to improve the quality of teaching and learning of mathematical modeling, particularly in the context of sustainability~\citep{lucas2023the,karjanto2023mathematical} 		

Although \cite{roe2018mathematics} have developed the theme of mathematics for sustainability excellently, that is, building progressive topics from measuring to deciding, we observed that the mathematical contents are rambling. In the beginning, a lot of algebra for solving physics-related problems were involved, then leap to system theory, graph theory, probability theory, and finally game theory. Calculus and differential equations were notably absent, but many concepts from economics appear intermittently, such as the market paradigm, Pareto efficiency, market failure, expected utility theory, and prospect theory.

In enhancing the teaching of mathematical modeling for sustainability, it is not only desirable to deviate from the textbook but also necessary to include other relevant mathematical topics and classroom examples that stimulate learning. For example, modeling of growth/decay has a direct implication with sustainability, and this topic requires both calculus and ordinary differential equations. Similarly, modeling interactions between species in an ecosystem reflects sustainability in ecological systems and promotes the coexistence of species in a way that benefits both the environment and human societies. Studying such interactions, such as the predator-prey system and population dynamics of species competition, requires a system of difference or differential equations, depending on whether the interest is in a discrete or continuous model.

It is essential to acknowledge that the diversity of mathematical topics within the context of sustainability extends beyond algebra, economics, and ecological modeling. To offer a comprehensive understanding of mathematical modeling's role in sustainability, we should explore additional domains. For instance, tackling resource allocation problems through optimization techniques is fundamental in addressing sustainability challenges, from efficient energy distribution to sustainable urban planning.

Additionally, the incorporation of real-world case studies and collaborative projects can enhance students' engagement and practical application of mathematical modeling for sustainability. By immersing learners in scenarios involving renewable energy adoption, waste reduction strategies, or urban resilience planning, they can grasp the tangible impact of mathematical modeling on real-world sustainability issues. Furthermore, fostering interdisciplinary collaboration among students from diverse academic backgrounds can replicate the holistic nature of sustainability problem-solving and prepare them for future roles in interdisciplinary teams dedicated to sustainable development.

In the journey to refine the teaching of mathematical modeling for sustainability, it is evident that an evolving and adaptable approach is necessary. By weaving together mathematical theory, interdisciplinary application, and practical examples, we aim to empower the next generation with the knowledge and skills required to contribute meaningfully to a more sustainable world.

\subsection*{Conflict of Interest}
The author declares no conflicts of interest.

\end{document}